\documentstyle{article}
\title{GG-functions and their relations to ${\cal A}$-hypergeometric functions}
\author{I.M.Gelfand and M.I.Graev
\thanks{The second author is supported by the Russian 
Foundation for Basic Research (grant 98--01--00798)}}

\begin{document}
\newtheorem{Def}{Definition}
\newtheorem{Th}{Theorem}
\newtheorem{Prop}{Proposition}
\setcounter{section}{-1}

\maketitle

\renewcommand{\theequation}{\thesection.\arabic{equation}}
\newcommand{\CC}{ {\bf C} }
\newcommand{\ZZ}{ {\bf Z} }
\newcommand{\RR}{ {\bf R} }
\newcommand{\om}{\omega}
\newcommand{\A}{ {\cal A} }

\section{Introduction}

In \cite{Gelf} I.M.Gelfand introduced the conception of hypergeometric
functions, associated with the Grassmanian $G_{k,n}$ of $k$-dimensional
subspaces in $\CC^n$. The class of these functions includes
many classical hypergeometric functions as particular cases.
Series of consequent papers of I.M.Gelfand
and his coathors (\cite{G1}--\cite{GVZ}) 
are devoted to the theory of these
functions. In particular, the paper \cite{G36} is devoted to 
hypergeometric functions associated with the Grassmanian $G_{3,6}\,$.

It should be mentioned that some well-known facts concerning classical
hypergeometric functions were simply interpreted in these investigations.
For example: 24 Kummer relations are known for the Gauss hypergeometric
function $F(a,b,c;x)$, which was proved to be connected with the
Grassmanian $G_{2,4}\,$. These relations arise from the natural action
of the permutation group $S_4$ at $G_{2,4}\,$. The Appel function $F_1$
is connected with the Grassmanian $G_{2,5}\,$. The existence of two integral
representations for $F_1$ by Euler integrals (one by single and one by double
integrals) arises from the isomorphism of Grassmanians $G_{2,5}$
and $G_{3,5}\,$.

Further development of the theory of hypergeometric functions is associated
with the paper of I.M.Gelfand, M.I.Graev and A.V.Zelevinsky \cite{GGZ}.
In this paper general hypergeometric systems of equations 
were defined and their holonomity was proved. These general
hypergeometric systems are called $\A$-hypergeometric systems
or GGZ-systems.

Any $\A$-hypergeometric system is defined by a set
$\A=\{\omega^1,\dots,\omega^N\}$
of vectors of ${\bf Z}^n$ that linearly generate $\CC^n$, and by a vector
$\beta\in\CC^n$. It consists of the following equations for functions
on $\CC^N$:
\equation\label{HS1}
\sum_{j=1}^N\left( a_j{\partial f \over\partial a_j}\right)
\om^j=f\cdot\beta
\endequation
\equation\label{HS2}
\prod_{\ell_j>0}\left({\partial \over\partial a_j}\right)^{\ell_j}f=
\prod_{\ell_j<0}\left({\partial\over\partial a_j}\right)^{-\ell_j}f
\endequation
for any $(\ell_1,\dots,\ell_N)\in\ZZ^N$ such that
$\sum_{j=1}^N\ell_j\omega^j=0\,$. Analytic solutions of this system
are called $\A$-hypergeometric functions.

Equation (\ref{HS1}) is the condition of homogeneity of the function
$f$ with respect to the following action of the complex torus $(\CC^*)^n$
at $\CC^N$:
$$
t:(a_1,\dots,a_n)\mapsto (t^{\om^1}a_1,\dots,t^{\om^N}a_N)\ ,
$$
where $t^{\omega^i}=t_1^{\omega^i_1}\cdot\dots\cdot t_n^{\om^i_n}\,$,
$t=(t_1,\dots,t_n)\in (\CC^*)^n\,$. By means of these conditions
one can reduce the system (\ref{HS1}),(\ref{HS2}) to a system of
$r=N-n$ equations for functions defined on the variety of orbits
of the torus $(\CC^*)^n$ in $\CC^N$, i.e., for functions of $r$ arguments.

All classical hypergeometric functions and all hypergeometric functions
associated with Grassmanians are $\A$-hypergeometric
functions for appropriate sets $\A$.

For example, the Gauss hypergeometric function is $\A$-hypergeometric
function associated with a set $\A$ of four vectors $\om^1$, $\om^2$,
$\om^3$, $\om^4$ of 3-dimensional space, connected by the single linear
relation $\om^1+\om^2-\om^3-\om^4=0\,$.

In this paper we describe the results of the papers \cite{GGf1}--\cite{GGf3},
where a new approach to the notion of general hypergeometric function was
suggested.

Let us illustrate this approach at the example of Gauss hypergeometric
function.

Consider the function
$$
f(a,b,c;x)=\sum_{n=0}^\infty
 {\Gamma(a+n)\Gamma(b+n)\over \Gamma(c+n) n!}x^n
$$
(here $\Gamma$ is the Euler Gamma-function). This function differs from
the Gauss function only by a constant multiplier. 
As a function of four arguments $a$, $b$, $c$, $x$,
this function satisfies the following relations:
$$
{d\over dx}f(a,b,c;x)=f(a+1,b+1,c+1;x)
$$
\equation\label{rels}
af(a,b,c;x)+xf(a+1,b+1,c+1;x)=f(a+1,b,c;x)
\endequation
$$
bf(a,b,c;x)+xf(a+1,b+1,c+1;x)=f(a,b+1,c;x)
$$
$$
cf(a,b,c;x)+xf(a+1,b+1,c+1;x)=f(a,b,c-1;x)
$$
(the last three relations are called Gauss relations).

The Gauss differential equation for the function $f$ 
\equation\label{Gauss}
x(1-x){d^2f\over dx^2}+\Bigl(c-(a+b+1)x\Bigr){df\over dx}-abf=0
\endequation
is a corollary
of the equations (\ref{rels}). Therefore for definition and investigation
of the function $f$ one can start not with a second-order differential
equation (\ref{Gauss}), but with the relations (\ref{rels}).

Such approach can be applied to any $\A$-hypergeometric system.
The system of differential equations can be replaced by a system of 
linear relations between a function $F\,$, its first-order partial 
derivatives, and its shifts with repect to parameters.  Namely,
instead of $\A$-hypergeometric system (\ref{HS1}),(\ref{HS2})), 
we introduce a system of equations in the space of functions $f(\beta,a)$ 
on $\CC^n\times\CC^N$. This system consists of equations (\ref{HS1}) 
and the following differential-difference equations:
\equation\label{difdif}
{\partial f(\beta,a)\over\partial a_j}=f(\beta-\om^j,a)\, ,\ j=1,\dots,N
\endequation
here the vectors $\om^i$ are not presumed to be integers.

System of equations (\ref{HS1}),(\ref{difdif}) is introduced and studied 
in \cite{GGf1}--\cite{GGf3}. There it is called GG-system associated
with $\A=\{\om^1,\dots,\om^N\}\,$. Its analytical with
respect to $\beta$ and $a$ solutions are called GG-functions 
associated with $\A\,$. Later, independently, and in a somewhat 
different form, these functions were introduced and studied in \cite{AI}. 
In \cite{AI} they are called quasi hypergeometric functions.

If vectors $\om^i$ are integer, then evidently (\ref{difdif})
implies (\ref{HS2}). In general case a GG-system cannot be reduced to
$\A$-hypergeometric systems, therefore solutions of GG-systems form
a wider class of functions than $\A$-hypergeometric functions.
However, the structure of these solutions (their
representations in the form of power series and in the form of integrals)
is proved to be similar to the structure of solutions of
$\A$-hypergeometric functions.

One of advantages of the new approach is that ``arguments" $a$ and 
``parameters" $\beta$ appear as arguments of same rights.  This allows 
us, in particular, to regard solutions of GG-systems 
that are distributions with respect to $\beta$ (see section 5).

In section 8 we give a general definition of GG-systems associated 
with an arbitrary complex Lie group. This definition includes 
the initial one as a particular case. 

\section{GG-systems and GG-functions}\setcounter{equation}{0}
\subsection{Definition of GG-system}
Let $V$ be an $n$-dimensional complex vector space and let
$$
\A = \{ \omega ^1, \dots , \omega ^N \}\, ,\ N\ge n
$$
be an arbitrary finite set of vectors from $V$ that linearly
generates $V\,$.
\begin{Def}\label{defGG}
GG-system, associated with the set $\A$, is the following
system of equations in the space of functions
$f(\beta,a)$ on $V\times\CC^N$:
\equation\label{GG1}
{\partial f(\beta,a)\over\partial a_j}=f(\beta-\omega^j,a)\, ,\ j=1,\dots,N
\endequation
\equation\label{GG2}
\sum_{j=1}^N a_j{\partial f\over\partial a_j}\cdot\omega^j=
  f\cdot\beta
\endequation
Solutions of this system in the class of analytical functions of $\beta$
and $a$ are called GG-functions associated with $\A\,$.
\end{Def}

\noindent{\bf Notes.}\ 1. From this definition it follows,
that GG-functions form a module over
the ring of functions $u(\beta)$, satisfying the periodicity conditions:
$$
u(\beta+\omega^j)=u(\beta)\, ,\ j=1,\dots,N\ .
$$

\noindent 2. By virtue of (\ref{GG1}), the equation (\ref{GG2}) can be
replaced by the following equation:
\equation\label{GG3}
\sum_{j=1}^N a_jf(\beta-\om^j,a)\cdot\om^j=f\cdot\beta
\endequation

\subsection{An equivalent definition of GG-system}
Let $L\subset\CC^n$ be an arbitrary fixed linear subspace and let
$L^\bot\subset(\CC^N)'$ be the ortogonal complement to $L$, i.e.,
the set of all $\nu\in(\CC^N)'$ such that $\langle\nu,\ell\rangle=0$
for all $\ell\in L\,$.
\begin{Def}\label{altdef}
GG-system, associated with the linear subspace $L\subset\CC^N$,
is the following system of equations in the space of functions
$F(\gamma,a)$, $\gamma,a\in\CC^N\,$:
\equation\label{alt1}
{\partial F(\gamma,a)\over\partial a_i}=F(\gamma-e_i,a)\, ,\ i=1,\dots,N\, ,
\endequation
where $e_1,\dots,e_N$ is the standart basis of $\CC^N$,
\equation\label{alt2}
F(\gamma+\ell,a)=F(\gamma,a)\quad\mbox{for any}\ \ell\in L
\endequation
\equation\label{alt3}
\sum_{i=1}^N\nu_ia_iF(\gamma-e_i,a)=\langle \nu,\gamma \rangle
F(\gamma,a)
\endequation
for any $\nu\in L^\bot$.
\end{Def}

Let us establish a connection between two definitions \ref{defGG}
and \ref{altdef}. To a set $\A=\{\om^1,\dots,\om^N\}$
of vectors of $V$ we correspond the following linear map:
$$
\pi:\CC^N\to V\ ,\quad\pi(e_i)=\om^i\ ,\quad i=1,\dots,N
$$
where $e_1,\dots,e_N$ is the standart basis of $\CC^N$.
Define $L=\ker\pi$. Conversly, if a linear subspace $L\subset\CC^N$
is given, then set $V=\CC^N\bigm/L$ and define $\om^i$, $i=1,\dots,N$ as
the image of $e_i$ under the natural projection $\CC^N\to V\,$. 

It is evident that if $f(\beta,a)$ is a function on $V\times\CC^N$ and
$F(\gamma,a)$ is its pull-back to $\CC^N\times\CC^N$, then the system
(\ref{GG1}),(\ref{GG2}) for $f$ is equivalent to the system
(\ref{alt1})--(\ref{alt3}) for $F\,$.

\subsection{Elementary and reducible cases}

{${\bf 1}^{\rm o}$}
\ There exists a unique, up to a coefficient, GG-function, associated
with $L=\CC^N$, namely $F(\gamma,a)=e^{a_1+\dots+a_N}$.

\noindent{${\bf 2}^{\rm o}$}
\ Any GG-function associated with the subspace $L=0$ has the form
$$F(\gamma,a)=u(\gamma)\prod_{j=1}^N{a_i^{\gamma_i}\over\Gamma(\gamma_i+1)}$$
where $u(\gamma)$ is an arbitrary periodic function with the period 1
with respect to every $\gamma_i\,$.

\noindent{${\bf 3}^{\rm o}$}
It can easily be checked that the GG-system on ${\bf C}^N\times{\bf C}^N$ 
can be reduced to a GG-system on a space of lower dimension if 
one of the following conditions is satisfied:
\begin{enumerate}
\item The subspace $L\subset\CC^N$ contains a nonzero coordinate subspace
\item The subspace $L$ is contained in a proper coordinate subspace
\item $L$ contains at least one vector of the form $e_i-e_j$, $i\ne j$,
where $\{e_i\}$ is the standart basis of $\CC^N$.
\end{enumerate}
In the sequel we may exclude these cases out of consideration.

In terms of the definition \ref{defGG} this amounts to the
following: first, the cases $n=0$ and $n=N$ are excluded; second,
vectors $\om^i$ are presumed to be non-zero, pairwise different, and
such that every $\om^i$ can be represented as a linear combination of others.

\subsection{Reduced GG-systems}
We reduce a GG-system, associated with a set $\A=\{\om^1,\dots,\om^N\}$
of vectors of $n$-dimensional space $V\,$, to a system of $N$ equations
for a function of $r=N-n$ arguments.

Fix an arbitrary solution $v(\beta,a)$ of the system
(\ref{GG2}) and an arbitrary basis $\{l^i=(l^i_1,\dots,l^i_n)|i=1,\dots,r\}$,
$r=N-n$ of the space $L\,$.

\noindent{\bf Lemma 1}\ \it Any solution $f(\beta,a)$ of the system
(\ref{GG2}) can be represented in the following form:
\equation\label{form}
f(\beta,a)=v(\beta,a)F(\beta,x)
\endequation
where $x=(a^{l^1},\dots,a^{l^r})\,$, $a^{l^j}=a_1^{l^j_1}\cdot\dots
\cdot a_N^{l^j_N}$.\rm
\smallskip

\noindent{\bf Lemma 2}\ \it Functions
$$\frac{\partial v(\beta,a)/\partial a_i}
{v(\beta-\om^i,a)}\quad\mbox{and}\quad
\frac{v(\beta,a)}{a_iv(\beta-\om^i,a)}\quad i=1,\dots,N$$
satisfy the equation (where $\varphi$ is any of these functions):
\equation\label{phieq}
\sum_{i=1}^N\left(a_i{\partial \varphi\over\partial a_i}\right)\om^i=0
\endequation
\rm\smallskip

\noindent{\bf Corollary}\ \it These functions can be represented in the form:
\equation\label{phipsi}
\frac{\partial v(\beta,a)/\partial a_i}{v(\beta-\om^i,a)}=\varphi_i(\beta,x)
\endequation
$$
\frac{v(\beta,a)}{a_iv(\beta-\om^i,a)}=\psi_i(\beta,x)
$$
where $x=(a^{l^1},\dots,a^{l^r})\,$.\rm
\smallskip

\begin{Th}
Any GG-function $f(\beta,a)$ associated with $\A=\{\om^1,\dots,\om^N\}$
can be represented in the form (\ref{form}), where $v$ is an arbitrary
fixed solution of the system (\ref{GG2}) and $F(\beta,x)$ is a solution
of the following system
\equation\label{varred}
F(\beta-\om^i,x)=\varphi_i(\beta,x)\cdot F(\beta,x)+
\sum_{j=1}^r l^j_i\psi_i(\beta,x)x_j{\partial F\over\partial x_j}(\beta,x)\ ,
\ i=1,\dots,N\ ;
\endequation
functions $\varphi_i$ and $\psi_i$ are defined by equations (\ref{phipsi}).

Conversely, if $F$ is a solution of the system (\ref{varred}),
then the function $f(\beta,a)\,$, defined by the equation (\ref{form}),
satisfies the GG-system (\ref{GG1}),(\ref{GG2}).
\end{Th}

\begin{Def}
We call the system (\ref{varred}) reduced GG-system associated
with $\A$, and its solutions in the class of
analytical functions of $\beta$ and $x$  
reduced GG-functions associated with $\A\,$.
\end{Def}

Let us emphasize, that reduced GG-system depends not only on $\A\,$,
but also on the choice of the solution $v$ of the system (\ref{GG2})
and of the basis $\{l^1,\dots,l^r\}$ in $L\,$.

We shall indicate a function $v$ and a basis of $L$ such that 
the equation (\ref{varred}) has the simplest form.

\begin{Def}
An $n$-element set $I\subset [1,N]$ is called a base if vectors
$\om^i\in V$, $i\in I$ are linearly independent (and hence form a basis
of $V$).
\end{Def}
Note that $I$ is a base if and only if $C^I\cap L=0$, where
$\CC^I\subset\CC^N$ is the coordinate subspace generated by the vectors
$e_i\,$, $i\in I\,$.

With any base $I$ we associate the linear map
$$\gamma:V\to\CC^N$$
defined by the equalities
$$\gamma(\omega^i)=e_i\ ,\ i\in I$$
Let us denote the $i$-th coordinate of $\gamma(\beta)$ by $\beta_i$ 
(so $\beta=\sum_{i\in I}\beta_i\om^i$). It is evident that the function
\equation\label{v}
v(\beta,a)=a^{\gamma(\beta)}=\prod_{i\in I}a_i^{\beta_i}
\endequation
satisfies the equation (\ref{GG2}). It is evident also that
the vectors
\equation\label{l}
l^j=e_j-\gamma(\om^j)\ ,\ j\in J=[1,N]\setminus I
\endequation
belong to the subspace $L$ and are linearly independent. Therefore
they form a basis of $L\,$.
\begin{Def}
The system (\ref{varred}), where the function $v$ and the basis
$\{l^j\}$ are defined by the equations (\ref{v}) and (\ref{l}), 
is called reduced GG-system associated with $\A$ and the base $I$.

Solutions of this system are called reduced GG-functions.
\end{Def}

\begin{Prop}
The reduced GG-system associated with $\A$ and $I$ has the following
form:
\equation\label{red1}
\beta_iF(\beta,x)+\sum_{j\in J} l_i^jx_j
{\partial F(\beta,x)\over\partial x_j}=
F(\beta-\om^i,x)\quad \mbox{for}\ i\in I
\endequation
where $l^j_i$ are coordinates of the vector $l^j\in\CC^N$, i.e.,
$l_i^j=-\gamma_i(\omega^j)\,$;
\equation\label{red2}
{\partial F(\beta,x)\over\partial x_j}=F(\beta-\omega^j,x)
\quad \mbox{for}\quad j\in J
\endequation
\end{Prop}

\section{Description of solutions of GG-systems}\setcounter{equation}{0}
\subsection{Solutions of a reduced GG-system that are regular in
a neighborhood of the point $x=0$}

For an arbitrary base $I\subset[1,N]$ we shall use the following
notation:
$$
\Gamma_I(\beta)=\prod_{i\in I}\Gamma(\beta_i)\, ;\ \beta\in V\, ,
$$
where $\Gamma$ is the Euler Gamma-function and $\beta_i$ are the coordinates
of the vector $\beta$ with respect to the basis $\{\omega^i\}_{i\in I}$
of $V\,$.

We shall write ${\bf 1}=(1,1,\dots,1)\in V\,$.

\begin{Th}\label{Th2}
Suppose $F(\beta,x)$ is a solution of the reduced GG-system 
(\ref{red1}),(\ref{red2}), associated with a set $\A$ and a base $I$. 
If $F(\beta,x)$ is regular in a neighborhood 
of the point $x=0$, then it has the following form:
\equation\label{GGsol}
F(\beta,x)=\sum_m
{u\left( \beta-\sum_{j\in J}m_j\om^j\right)
\over\Gamma_I\left( \beta-\sum_{j\in J}m_j\om^j+{\bf 1}\right)}
\cdot{x^m\over m!}\ ,
\endequation
where $m=(m_j)_{j\in J}\,$, $m_j\in\ZZ_+\,$, 
${\displaystyle {x^m\over m!}=\prod_{j\in J}{x_j^{m_j}\over m_j!}}\,$, and
$u(\beta)$ is a function on $V$ that satisfies the periodicity condition:
\equation\label{period}
u(\beta-\om^i)=u(\beta)\quad\mbox{for any}\ i\in I\,.
\endequation
                   
Conversely, if $F$ is a formal series of the form (\ref{GGsol}), where
$u(\beta)$ is an arbitrary function satisfying the periodicity
condition (\ref{period}), then $F$ formally satisfies the equations of
reduced GG-system (\ref{red1}),(\ref{red2}).
\end{Th}

\noindent{\bf Corollary.}\it
\ Any solution of reduced GG-system (\ref{red1}),
(\ref{red2}), regular in a neighborhood of $x=0$, is determined
uniquely by the initial term $u(\beta)/\Gamma_I(\beta+{\bf 1})$
of its power series expansion in powers of $x\,$.\rm

\begin{Def}
We denote formal power series (\ref{GGsol}), where 
$$
u(\beta)=\exp(2\pi i\langle k,\beta\rangle )\ ,
\ \langle k,\beta\rangle=\sum_{i\in I} k_i\beta_i\ k\in \ZZ^I\, ,
$$
by $F_{I,k}\,$. 
We call them reduced 
GG-series, associated with $\A$ and a base $I\subset [1,N]$. 
\end{Def}
It follows from the definition that
$$
F_{I,0}(\beta, x)=
\sum_m {x^m\over\Gamma_I\left(\beta-\sum_{j\in J}m_j\om^j+{\bf 1}\right) m!}
$$
$$
F_{I,k}(\beta,x)=\exp(2\pi i\langle k,\beta\rangle)
F_{I,0}\left(\beta,e^{-2\pi i\langle k,\om^J\rangle}x\right)\ ,\ k\in{\bf Z}^I
$$
where $e^{-2\pi i\langle k,\om^J\rangle}x$ is the vector with coordinates
$e^{-2\pi i\langle k,\om^j\rangle}x_j\,$, $j\in J\,$.

By $l^j_i\,$, $i\in I$ denote coordinates of the vector
$\om^j$, $j\in J$ with respect to the basis $\{\om^i\}_{i\in I}\,$.

\begin{Prop}
If
\equation\label{condi}
{\rm Re}\,\sum_{i\in I}l^j_i\ge -1\ \mbox{for all}\ j\in J\ ,
\endequation
then the series $F_{I,k}$ converge in a neighborhood of $x=0$ and thus
define reduced GG-functions, which are regular in this neighborhood .
\end{Prop}

\noindent{\bf Remark} The series (\ref{GGsol}) can be represent 
also in the form 
$$
F(\beta,x)=\sum_m
u\left( \beta-\sum_{j\in J}m_j\om^j\right)
{\Gamma_{I_1}\left( -\beta+\sum_{j\in J}m_j\om^j\right)
\over\Gamma_{I_2}\left( \beta-\sum_{j\in J}m_j\om^j+{\bf 1}\right)}
\cdot{x^m\over m!}\ ,
$$
where $I=I_1\sqcup I_2$ is an arbitrary fixed partition of $I$, 
and $u$ satisfies the following conditions:
$$u(\beta-\omega^i)=-u(\beta)\ \mbox{if}\ i\in I_1\,,
\ u(\beta-\omega^i)=u(\beta)\ \mbox{if}\ i\in I_2
$$
\subsection{Solutions of GG-system (\ref{GG1}),(\ref{GG2})
that are regular in a neighborhood of coordinate subspaces of
$\CC^n$.}
Let $I\subset[1,N]$ be an arbitrary base and $\CC^I\subset\CC^N$
be the coordinate subspace generated by vectors $e_i\,$, $i\in I\,$.
We say that $a\in\CC^I$ is a generic point if $a_i\ne 0$
for all $i\in I\,$.

Theorem (\ref{Th2}) implies 
\begin{Th}
Any solution of GG-system (\ref{GG1}),(\ref{GG2}) that is regular
in a neighborhood of a generic point of $C^I$, where $I\subset[1,N]$ is
an arbitrary base, has the following form:
\equation\label{aGG}
f(\beta,a)=\prod_{i\in I}a_i^{\beta_i}\sum_m
\frac{u\left(\beta-\sum_{j\in J}m_j\om^j\right)}
{\Gamma_I\left(\beta-\sum_{j\in J}m_j\om^j+1\right)}
\prod_{j\in J}{y_j^{m_j}\over m_j!}
\endequation
where $u$ satisfies (\ref{period}); $J=[1,N]\setminus I$, 
$y_j=a_j\prod\limits_{i\in I}a_i^{l_i^j}$; 
$\beta_i$ and $l_i^j$ are coordinates of 
$\beta$ and $-\om^j$, respectively, 
with respect to the basis $\{\om^i\}_{i\in I}\,$.

Conversely, if $f$ is a formal series of the form (\ref{aGG}), where
$u(\beta)$ satisfies (\ref{period}), 
then $f$ formally satisfies the equations of GG-system (\ref{GG1}),(\ref{GG2}).
\end{Th}
\begin{Def}
By $f_{I,k}(\beta,a)\,$, $k\in\ZZ^I$, we denote 
the formal series (\ref{aGG}), where $u(\beta)=
\exp(2\pi i\langle k,\beta\rangle)\,$.
We say that these series are GG-series, associated with $\A$ and $I$,
\end{Def}

It follows from the definition that GG-series and reduced GG-series
are linked by the relation:
$$
f_{I,k}(\beta,a)=a^\beta F_{I,k}(\beta,x)
$$
where $a^\beta=\prod_{i\in I}a_i^{\beta_i}$ and $x\in\CC^J$ is the vector
with coordinates $x_j=a_j\prod_{i\in I} a_i^{l_i^j}$.

\section{Relation to $\A$-hypergeometric functions}\setcounter{equation}{0}
Consider the case when all vectors $\om^1,\dots,\om^N$ that
form the set $\A$ have integer coordinates with respect to some
fixed basis of $V\,$. Then with the set $\A$ one can 
associate both GG-system (\ref{GG1}),(\ref{GG2}) and
$\A$-hypergeometric system (\ref{HS1}),(\ref{HS2}).

We noted already that for any fixed $\beta$ a solution of the GG-system 
(\ref{GG1}),(\ref{GG2}) satisfies the $\A$-hypergeometric system (\ref{HS1}),
(\ref{HS2}).

\begin{Th}\label{set}
There exists a finite set $\{f_i(\beta,a)\}$ of GG-functions,
associated with $\A$ and defined on $V\times U\,$, where
$U\subset\CC^N$ is some domain, such that
\begin{enumerate}
\item Any solution of GG-system (\ref{GG1}),(\ref{GG2}) on $V\times U$
can be uniquely represented in the form:
$$f(\beta,a)=\sum_i u_i(\beta) f_i(\beta,a)$$
where functions $u_i(\beta)$ satisfy the periodicity condition:
$$u_i(\beta+\om^j)=u_i(\beta)\quad \mbox{for any}\quad j=1,\dots,N$$
\item For any fixed generic $\beta\in V$ the functions $f_i(\beta,a)$
form a basis in the space of solutions of $\A$-hypergeometric system
(\ref{HS1}),(\ref{HS2}) on $U\,$.
\end{enumerate}
\end{Th}

Let us present a construction of such set of functions $f_i\,$.

Let $L\subset\CC^N$ be the subspace of all vectors
$(\ell_1,\dots,\ell_N)\,$, satisfying the condition
$\sum\limits_{j=1}^N \ell_j\om^j=0\,$, and $\Lambda\subset(\ZZ^N)'$
be the lattice of integer vectors orthogonal to $L\,$.
By $\Lambda^I\subset (\ZZ^I)'$ denote the projection
of $\Lambda$ to $(\ZZ^I)'$, where $\ZZ^I$ is the lattice generated
by the vectors $e_i\,$, $i\in I\,$,                               . 
In the sequel, we denote
${\cal L}^I=(\ZZ^I)'\bigm/\Lambda^I\,$. For any base $I$ the set
${\cal L}^I$ is finite.

Consider the GG-series $f_{I,k}(\beta,a)$ from section 2.

Let ${\cal J}$ be a subset of the set of all bases. 
By ${\cal F}_{\cal J}$ denote
the set of GG-series $f_{I,k}\,$, where
$I$ runs through bases from ${\cal J}\,$, and $k$ for every base $I$
runs through one representative of every class
$\kappa\in{\cal L}=(\ZZ^I)'/\Lambda^I$.

\begin{Th}
There exist a subset ${\cal J}$ of bases and a domain
$U\subset\CC^N$ such that
\begin{enumerate}
\item all GG-series $f_{I,k}(\beta,a)\in{\cal F}_{\cal J}$
converge in $U$
\item the set ${\cal F}_{\cal J}$ of GG-functions satisfies
the conditions of the theorem \ref{set}.
\end{enumerate}
\end{Th}

\section{Integral representations of GG-functions}\setcounter{equation}{0}
\subsection{Formal solutions}
Consider a GG-system associated with a set ${\cal A}=
\{\om^1,\dots,\om^N\}$ of vectors of a space $V$ of dimension $n\,$.
Fix an arbitrary basis of $V$ and for $\beta\in V$
denote coordinates of $\beta$ with respect
to this basis by $\beta_1,\dots,\beta_n$. 
Also denote: $T=V\setminus\Sigma\,$, where
$\Sigma$ is the union of all coordinate subspaces of $V\,$.
\begin{Prop}
Integrals of the following form formally satisfy the GG-system
(\ref{GG1}),(\ref{GG2}):
\equation\label{GGint}
f(\beta,a)=\int\limits_C\exp\left(\sum_{j=1}^N a_jt^{\om^j}\right)
t^{-\beta}\prod_{i=1}^n{dt_i\over t_i}\ ;
\endequation
here $t^{\om^j}=t_i^{\om^j_1}\dots t_n^{\om^j_n}$,
$t^{-\beta}=t_1^{-\beta_1}\dots t_n^{-\beta_n}$ and $C\subset T$
is an arbitrary cycle of real dimension $n\,$.
\end{Prop}

By the change of variables $t_i=e^{s_i}$ the integral (\ref{GGint})
can be represented in the following form:
\equation\label{int3}
f(\beta,a)=\int\limits_\Gamma\exp
\left(\sum_{j=1}^N a_je^{\langle\om^j,s\rangle}-\langle\beta,s\rangle\right)
ds_1\dots ds_n\ ,
\endequation
where $\Gamma$ is an $n$-dimensional cycle in the space $V'$
dual to $V\,$.

If the mentioned basis is the set $\{\om^i\}_{i\in I}\,$,
where $I\subset[1,N]$ is a base, then the integral (\ref{GGint})
has the following form:
\equation\label{int2}
f(\beta,a)=\int\limits_C\exp
\left(\sum_{i\in I}a_it^i+\sum_{j\notin I} a_jt^{\ell^j}\right)
t^{-\beta}\prod_{i\in I}{dt_i\over t_i}\ ,
\endequation
where $t^{\ell^j}=\prod_{i\in I} t_i^{\ell_i^j}$ and $\ell^j_i$
are coordinates of $\om^j$, $j\notin I$ with respect to the
choosen basis.
\begin{Prop}
For any base $I\subset [1,N]$ the following integrals formally satisfy
the reduced GG-system (\ref{red1}),(\ref{red2})
associated with the base $I$:
\equation\label{redint}
F(\beta,x)=\int\limits_C\exp
\left(\sum_{i\in I} t_i+\sum_{j\notin I}x_jt^{\ell^j}\right)
t^{-\beta}\prod_{i\in I}{dt_i\over t_i}\ ;
\endequation
here $C\subset T$ is an arbitrary cycle of real dimension $n\,$.
\end{Prop}

\noindent{\bf Remark.}\ If $\sum\limits_{i\in I}l_i^j=1$ for every 
$j\notin I$ then the following integrals formally satisfy also 
system (\ref{red1}),(\ref{red2}): 
\equation
F(\beta,x)=\Gamma^{-1}(\beta_I+1)\int\limits_C
\left(1+\sum_{j\notin I}x_js^{l^j}\right)^{\beta_I}
\prod_{i\in I}s_i^{-\beta_i-1}ds\, ,
\endequation
where $\beta_I=\sum\limits_{i\in I}\beta_i$ and $C$ is an arbitrary cicle 
of real dimension $n-1$ in the hyperplane $\sum\limits_{i\in I}s_i=1\,$. 

Similarly to (\ref{GGint}), the integral (\ref{redint}) can be
represented in the form:
\equation\label{redint2}
F(\beta,x)=\int\limits_\Gamma\exp
\left(\sum_{i\in I} e^{s_i}+\sum_{j\notin I}x_je^{\langle\om^j,s\rangle}
-\langle\beta,s\rangle\right)ds_1\dots ds_n
\endequation

The problem is to find a family of cycles such that these integrals
have nonempty domain of convergency. One may solve this problem
only for the integrals (\ref{redint}).

\subsection{Case $n=1$.}
In the case $n=1$ the set ${\cal A}$ is a set of nonegative
numbers $\om^1,\dots,\om^N\in\CC\,$, and bases are one-element subsets
of $[1,N]\,$. Integral (\ref{redint}) associated with an arbitrary
base $I=\{i\}$ has the form:
\equation\label{int1}
F(\beta,x)=\int\limits_C\exp\left(t+\sum_{j\ne i}x_jt^{\om_j/\om_i}\right)
t^{-\beta-1}dt
\endequation
\begin{Prop}\label{1dim}
Let $C\in\CC$ be a loop, going out of $-\infty$ by
the real negative axis, around 0,
and again to $-\infty$ by the negative axis.
If ${\rm Re}(\om_j/\om_i)\le 1\,$, then the integral (\ref{int1})
converge in a neighborhood of $x=0$ and coincides with the function:
$$
F_0(\beta,x)=2\pi i\sum_m
\left(
\Gamma^{-1}\Bigl(\beta-\sum_{j\ne i}{\om_j\over\om_i}m_j+1\Bigr)
\prod_{j\ne i}{x_j^{m_j}\over m_j!}
\right)
$$
\end{Prop}
\subsection{Case $n>1\,$.}
Let us regard the integral (\ref{redint}) connected with an arbitrary
base $I\subset[1,N]$.
\begin{Prop}
If $\sum_{i\in I}\mathop{\rm Re}\ell^j_i\le 1$ for all
$j\in[1,N]\setminus I$ then there exists a cycle such that the integral
(\ref{redint}) converges in a neighborhood of the point $x=0\,$.
\end{Prop}
Let us present a description of the cycle $C\,$. Let us regard the
change of variables $t=\rho u\,$, where $\rho\in\CC\setminus\{0\}$
and $u$ is a point of the plane $U=\{u\in T|\sum u_i=1\}\,$.
By means of this change the integral (\ref{redint}) can be reduced to the
integral of the following form:
\begin{eqnarray}
\lefteqn{F(\beta,x)=\int\limits_C\exp\left(
\rho+\sum_{j\notin I}\left(x_j\rho^{r_j}\prod_{i\in I}u_i^{\ell^j_i}\right)
\right)\times}&&\label{rhoint}\\
&&\times\rho^{-\sum\beta_i-1}d\rho\,\prod_{i\in I}u_i^{-\beta_i-1}du
\nonumber
\end{eqnarray}
where $r_i=\sum\limits_{i\in I}\ell^j_i$ and $du$ is a holomorphic
volume form on $U\,$. We shall describe the cycle $C$ in coordinates
$\rho$ and $u\,$.

To give a precise meaning to the integral (\ref{rhoint}) let us
replace the multifunction on $U$ under the integral by a function
on a ramified covering $\tilde U\cong\CC^{n-1}$ over $U\,$.
For this end let us introduce an infinite covering $\tilde T\cong\CC^I$
over $T$ such that the prototype of a point $t\in T$ consists of
points $s\in\CC^I$ defined by the equality $e^{s_i}=t_i\,$, $i\in I\,$.
We define $\tilde U$ as the prototype of $U$ under this covering.
We shall interpretate the multiform from integral (\ref{rhoint})                      
as single-valued form on $\tilde U\,$.

Let us introduce a cycle $\tilde\Gamma\subset\tilde U\,$, which is
a multidimensional analog of double loop. The cycle $\tilde\Gamma$
is associated with the simplex
$$\Delta=\{u\in U\cap\RR^I|u_i>0\}\ ,$$
its construction see in \cite{Loop1},\cite{Loop2}. The desired cycle
$C\subset(\CC\setminus\{0\})^n\times\tilde U$ is equal to
$C_1\times\tilde\Gamma\,$,
where $C_1\subset\CC\setminus\{0\}$ is a loop around 0 along the negative
axis (the same as $C$ from the proposition \ref{1dim}).

\subsection{Integral representation of GG-functions associated
with real vectors.}
With any base $I\subset[1,N]$ we associate the following $n$-dimensional
cycles $\Gamma_k\subset\CC^I$, $k\in\ZZ^I$:
$$
\Gamma_k=\{s\in\CC^I|\mathop{\rm Im} s_i=(2k_i+1)\pi,\,i\in I\}
$$
\begin{Prop}
If all vectors $\om^1,\dots,\om^N$ are real (with respect to some
basis in $V$), then the integral (\ref{redint2}) over the cycle $\Gamma_k$
converges in the domain
$$\mathop{\rm Re}\beta_i>0\,,\ i\in I\, ,$$
$$
\left|\mathop{\rm arg}\left(x_je^{\pi i\langle \om^j,2k+{\bf 1}\rangle}
\right)\right|>{\pi\over 2}\ ,\ j\in [1,N]\setminus I
$$
and thus is a reduced GG-function.
\end{Prop}
We denote this function by $F_k(\beta,x)$.
\smallskip

\noindent{\bf Note.}\ The function $F_K$ can be extended by analiticity
to all values of $\beta\,$.

\section{Special cases}
\subsection{Case $N=n+1$}. In this case the reduced GG-functions are 
functions $F(\beta,x)$ on ${\bf C}^N\times{\bf C}\,$. The reduced 
GG-system has the following form:
\equation\label{scone1}
\beta_iF(\beta,x)+\ell_ix{dF(\beta,x)\over dx}=
F(\beta-e_i,x)\ ,\ i=1,\dots,n
\endequation
\equation\label{scone2}
{dF(\beta,x)\over dx}=F(\beta+\ell,x)
\endequation
where $\ell=(\ell_1,\dots,\ell_N)$ is an arbitrary fixed vector in ${\bf C}^N$. 
Every solution of (\ref{scone1}),(\ref{scone2}) that is regular in 
a neighborhood of $x=0$ has the following form
$$
F(\beta,x)=\sum_{m=0}^\infty u(\beta+m\ell)c(m){x^m\over m!}\ ,
$$
where $c(m)=\left(\prod\limits_{i=1}^n\Gamma(\beta_i+m\ell_i+1)\right)^{-1}$, 
and $u$ is an arbitrary periodical function on ${\bf C}^n$ with period 1 
with respect to every $\beta_i\,$. 

Examples. If $n=1$ then 
$$
F(\beta,x)=\sum_{m=0}^\infty {u(\beta+m\ell)\over\Gamma(\beta+m\ell+1)}
{x^m\over m!}\ ,\ \beta,\ell\in{\bf C}
$$ 

If $n=2$ then 
$$
F(\beta_1,\beta_2,x)=\sum_{m=0}^\infty 
{u(\beta_1+m\ell_2,\beta_2+m\ell_2)\over
\Gamma(\beta_1+m\ell_1+1)\Gamma(\beta_2+m\ell_2+1)}
{x^m\over m!}\ .
$$ 

The integral representation of the GG-functions has 
the form
$$
F(\beta,x)=\int\limits_C\exp(t_1+\dots+t_n+t_1^{-\ell_1}\dots t_n^{-\ell_n}x)
\prod_{i=1}^nt_i^{-\beta_i-1}dt_i
$$
\subsection{Case $n=1$.} In this case the reduced GG-functions are 
functions $F(\beta,x_1,\dots,x_r)$ on ${\bf C}\times{\bf C}^r$. 
The reduced GG-system has the form:
\equation\label{sctwo1}
\beta F(\beta,x)+\sum_{j=1}^r\ell^jx_j{\partial F(\beta,x)\over \partial x_j}=
F(\beta-1,x)
\endequation
\equation\label{sctwo2}
{\partial F(\beta,x)\over \partial x_j}=F(\beta+\ell^j,x)
\endequation
where $\ell=(\ell^1,\dots,\ell^r)$ is an arbitrary fixed vector in ${\bf C}^r$. 
Every solution of (\ref{sctwo1}),(\ref{sctwo2}) that is regular in 
a neighborhood of $x=0$ has the following form
$$
F(\beta,x)=\sum_m \left({u(\beta+\sum\limits_{j=1}^rm_j\ell^j)\over
\Gamma(\beta+\sum\limits_{j=1}^rm_j\ell^j+1)}
\prod_{j=1}^r{x_j^{m_j}\over m_j!}\right)\ ,
$$
where $u$ is an arbitrary periodical function on ${\bf C}$ with period 1.

The integral representation of the GG-functions has 
the form
$$
F(\beta,x)=\int\limits_\Gamma
\exp(t+t^{-\ell^1}x_1+\dots +t^{-\ell^r}x_r)t^{-\beta-1}dt
$$
This integral converges for every $x$ in case when $\mathop{\rm Re}\beta>0\,$,
$\mathop{\rm Re}\ell^j>-1\,$, $j=1,\dots,r$ and 
$\Gamma$ is a loop going from $-\infty-i0$ to $-\infty+i0$ around the point 0. 

\section{GG-distributions}\setcounter{equation}{0}
\subsection{Definition of GG-distributions}
Let us denote the space of compactly
supported $C^\infty$--functions on $\RR^n$ by $K_n$, 
and the space of Fourier transforms
of functions from $K_n$ by $Z_n$, with natural topologies. 
According to Paley-Wiener theorem the space $Z_n$ consists of entire 
analytical functions $F(z_1,\dots,z_n)$, satisfying the estimation:
$$
|z_1^{q_1}\dots z_n^{q_n}F(z_1,\dots,z_n)|\le
C_qe^{a_1|\mathop{\rm Im}z_1|+\dots+a_n|\mathop{\rm Im}z_n|}
$$
for any $q_1,\dots,q_n=0,1,\dots$.

We say that elements of dual spaces $K'_n$ and $Z'_n$,
i.e., continious linear functionals on $K_n$
and $Z_n$, are {\em distributions} on $\RR^n$ and $\CC^n$, respectively.
In particular, any continuous function 
$f(\xi)$ on $\RR^n$ defines a continuous linear functional on $K_n$
by the following formula:
$$
(f,F)=\int\limits_{\RR^n}f(\xi)F(\xi)d\xi
$$
thus $f$ can be regarded as a distribution on $\RR^n$.

Suppose $f$ is a distribution on $\RR^n$. 
The Fourier transform of $f$ is the linear functonal $\tilde f$ on $Z_n$
such that for any function
$F\in K_n$ and its Fourier transform $\tilde F\in Z_n$
the following equality holds:
$$
\Bigl(\tilde{f}, \tilde{F}(z)\Bigr)=(2\pi)^n\Bigl(f,F(-\xi)\Bigr)\ .
$$
In particular, the Fourier transform of the function
$e^{i\langle a,\xi\rangle}$, where $a=(a_1,\dots,a_n)$,
$\langle a,\xi\rangle=\sum a_i\xi_i\,$, is the distribution
$(2\pi)^n\delta(z-a)\in Z'_n\,$, defined by the equality
$$\Bigl(\delta(z-a),\tilde{F}(z)\Bigr)=\tilde{F}(a)$$
for any $\tilde{F}\in Z_n\,$.

We denote the space of analytical functions on $\CC^n$
with values in $Z'_n$ by $H=H_{n,r}\,$. We denote elements of $H$ by
$\varphi(\beta,x)\,$, $\beta\in\CC^n$, $x\in\CC^r$.

Let us regard the reduced GG-system (\ref{red1}),(\ref{red2}) associated
with a set $\A=\{\om^1,\dots,\om^N\}$ of vectors of an $n$-dimensional
linear space $V$ and a base $I\subset[1,N]\,$.
Without loss of generality we may assume that $I=[1,n]$ and denote
$\om^i=e_i$, $i=1,\dots,n$ and $\om^{n+q}=\ell^q=\sum_{i=1}^n\ell_i^qe_i$
for $q=1,\dots,r\,$, $r=N-n\,$. Then the equations of GG-system
take the following form:
\equation\label{dis1}
{\partial\varphi(\beta,x)\over\partial x_q}=\varphi(\beta-e_q,x)\, ,
\ q=1,\dots,r
\endequation
\equation\label{dis2}
\varphi(\beta-e_p,x)+\sum_{q=1}^r\ell_p^q\varphi(\beta-\ell^q,x)=
\beta_p\varphi(\beta,x)\,,\ p=1,\dots,n
\endequation
where $\{e_i\}$ is the standart basis of $\CC^n$. It is evident, that
these equations have sence also for elements of the space $H\,$.
\begin{Def}
Elements of $H$ that satisfy (\ref{dis1}),(\ref{dis2})
are called GG-distributions associated with $\A$ and $I\,$.
\end{Def}
\subsection{Description of GG-distributions}
\begin{Th}\label{Dth1}
For any GG-system (\ref{dis2}),(\ref{dis2}) there exists a unique up 
to a constant multiplier GG-distribution $f\,$. This distribution is 
the Fourier transform with respect to $\xi=(\xi_1,\dots,\xi_n)$ 
of the following function:
$$
F(\xi,x)=\exp\left(\sum_{p=1}^n e^{i\xi_p}+\sum_{q=1}^r 
x_qe^{\langle \ell^q,\xi\rangle}\right)\,,
$$
where $\langle \ell^q,\xi\rangle=\ell^q_1\xi_1+\dots+\ell^q_n\xi_n\,$.

Thus formally $f$ is defined by the integral
$$
f(\beta,x)=\int_{\RR^n}\exp\left(e^{i\xi_p}+\sum_{q=1}^r 
x_qe^{\langle \ell^q,\xi\rangle}-i\langle\beta,\xi\rangle\right)
d\xi_1\dots d\xi_n
$$
\end{Th}

\subsection{Representation of the GG-distribution as series}
\begin{Th}
The GG-distribution $f(\beta,x)$ can be represented as the
following power series with respect to $x\,$:
$$
f(\beta,x)=\sum_{m\in\ZZ_+^k}\left(
c_m(\beta)\prod_{q=1}^k{x^{m_q}\over m_q!}\right)\ ,
$$
where
$$
c_m(\beta)=\sum_{r\in\ZZ^n_+}\left(\prod_{p=1}^n r_p!\right)^{-1}
\delta\Bigl(\beta_1-\sum_{q=1}^k\ell_1^qm_q-r_1,\dots,
\beta_n-\sum_{q=1}^k\ell_n^qm_q-r_n\Bigr)\ ,
$$
$\delta(\beta_1,\dots,\beta_n)$ is the delta-function on $\CC^n$.
\end{Th}

In other words, for any function $\varphi\in Z_n$ we have
\equation\label{Dser}
(f,\varphi)=\sum_{m\in\ZZ_+^k}(c_m,\varphi){x^m\over m!}\ ,
\endequation
where ${\displaystyle {x^m\over m!}=\prod_{q=1}^k{x_q^{m_q}\over m_q!}}\,$,
$$
(c_m,\varphi)=\sum_{r\in\ZZ_+^n}(r!)^{-1}
\varphi\Bigl(\sum_{q=1}^k\ell_1^qm_q+r_1,\dots,\sum_{q=1}^k\ell_n^qm_q+r_n
\Bigr)\ .
$$
\begin{Prop}
Series (\ref{Dser}) converges for any $\varphi\in Z_n$ and $x\in\CC\,$.
\end{Prop}

The distributions $c_m\in Z_n^*$ are analytical functionals, i.e.,
$(c_m,\varphi)$ can be represented in the form:
$$
(c_m,\varphi)=\int_{\Gamma_m}F_m(\beta)\varphi(\beta)\,d\beta
$$
where $F_m$ is a function and $\Gamma_m\subset\CC^n$ is a surface
of real dimension $n\,$. In the case $k=1$ we have:
$$
F_m(\beta)=\prod_{p=1}^n e^{\pi i(\ell_pm-\beta_p)}
\Gamma(\ell_pm-\beta_p)
$$
and $\Gamma_m=\gamma_m^1\times\dots\times\gamma^n_m\,$,
where $\gamma_m^p\subset\CC$ is a contour that goes around poles
of the function $\Gamma(\ell_pm-\beta_p)\,$.

\subsection{Connections with series of hypergeometric type}
There exists a formal operation $\varphi\mapsto f$ that refers
a GG-distribution $f$ to any GG-function $\varphi$ if $\varphi$
is defined by a series of hypergeometric type. This operation
is replacing the multipliers $\Gamma(u)$ and $1/\Gamma(u+1)$
in the coefficients of the series by,
respectively, the distributions $\gamma^+(u)$ and $\gamma^-(u)\,$,
where
$$
\gamma^+(s)=\sum_{m=0}^\infty{(-1)^m\over m!}\delta(s+m)\,,
\quad
\gamma^-(s)=\sum_{m=0}^\infty{1\over m!}\delta(s-m)\,;
$$
$\delta(s)$ is the delta-function on $\CC\,$.

The functions $\gamma^+(s)$ and $\gamma^-(s)$ satisfy the following
functional relations:
$$
\gamma^+(s+1)=s\gamma^+(s)\,,\quad\gamma^-(s-1)=s\gamma^-(s)\ ,
$$
i.e., the same relations as $\Gamma(s)$ and $1/\Gamma(s+1)$
respectively.

By theorem \ref{Dth1}, the images under this operation of all
hypergeometric-type series
that satisfy system (\ref{GG1}),(\ref{GG2}) differ only by
constant multipliers.

\section{Resonance GG-systems}\setcounter{equation}{0}
     
There exist GG-systems such that their solutions satisfy
some additional differential equations under certain
relations between parameters and certain regularity conditions.
We call such GG-systems {\em resonance}. 

\subsection{Resonance sets}
Let $\A=\{\om\}$ be a finite set of non-zero vectors linearly
generating a space $V\,$. Let us introduce the following
notation for an arbitrary vector $v\in V\,$:
$$
\A_v=\Bigl\{\om\in\A\bigm|\om+v\in\A\cup\{0\}\Bigr\}\,,\quad
{\cal B}_v=\A\setminus\A_v\ ,
$$
$L_v\subset V$ is the linear subspace generated by the vectors
$\om\in{\cal B}_v\,$. In particular, $\A_0=\A\,$, ${\cal B}_0=\emptyset\,$,
$L_0=0\,$.

\begin{Def}
A vector $v\in V$ is called consistent with the set $\A$ if 
$L_v$ is a proper subspace of $V\,$.
\end{Def}

In particular, the vector $v=0$ is consistent with any $\A\,$. 
Obviously, if a vector $v\ne 0$ is consistent with $\A\,$, 
then either $v=\om'-\om$ or $v=-\om\,$, where $\om,\om'\in\A\,$. 
Thus vectors in general position are not consistent with $\A\,$. 

Example: $\A=\{\om^1,\dots,\om^{n+1}\}\,$, where $\om^i$ are vectors 
from $\CC^n\,$, and $\om^1,\dots,\om^n$ are linearly independent. 
In this case all vectors $\om^{n+1}-\om^i$, $i=1,\dots,i$ as well 
as the vector $-\om^{n+1}$ are consistent with $\A\,$. 

It is easily proved that for any consistent with $\A$ vector $v\ne 0$
we have: $\mathop{\rm codim}L_v=1\,$, $v\notin\A\,$, $v\notin L_v\,$, 
and $A_v\cap L_v=\emptyset\,$. This implies that any set $\A$ of $N$ vectors 
of $V$ such that a fixed nonzero vector $v$ is consistent with $\A\,$, 
has the following form:
\equation\label{conc}
\A=\bigl\{\om^i-jv\bigm| j=0,1,\dots,k_i-1;i=1,\dots,r\bigr\}
\cup\bigl\{-jv\bigm|j=1,\dots,k_0-1\bigr\}\ ,
\endequation
where $\om^1,\dots,\om^r$, $r<N$ are arbitrary pairwise different vectors 
such that the space, generated by them, has codimension 1 and does not 
contain $v\,$, and $k_0,\dots,k_r$ are arbitrary natural numbers 
such that $\sum\limits_{i=1}^r k_i=N+1$ (if $k_0=1\,$, then the second set 
in (\ref{conc}) is assumed to be empty).

\subsection{Resonance GG-systems}
Let us regard a GG-system on $V\times\CC^N$, associated with a set 
$\A=\{\om\}$ of $N$ nonzero vectors of $V$, linearly generating $V\,$.
For convenience we shall enumerate coordinates in $\CC^N$ not by numbers, 
but by elements of $\A\,$, i.e., instead of $a_i\,$, $i=1,\dots,N$ 
we shall write $a_\om\,$, $\om\in\A\,$. In this notation the GG-system 
has the following form:
\equation\label{GGA1}
{\partial f(\beta,a)\over\partial a_\om}=f(\beta-\omega,a)\, ,\ \om\in\A
\endequation
\equation\label{GGA2}
\sum_{\om\in\A} a_\om{\partial f\over\partial a_\om}\cdot\om=
  f(\beta,a)\cdot\beta
\endequation

\begin{Def}
The GG-system {\rm (\ref{GGA1}),(\ref{GGA2})} is called resonance
if there exists at least one nonzero consistent 
with $\A$ vector $v\in V\,$. 
\end{Def}
\begin{Def}
Suppose $v\ne 0$ is an arbitrary consistent with $A$ vector. 
Let us associate with $v$ the following hyperplane:
$$
\Pi_v=L_v+v
$$
\end{Def}

\begin{Th}
Suppose GG-system {\rm (\ref{GGA1}),(\ref{GGA2})} is resonance, 
$v\in V$ is a nonzero consistent with $\A$ vector, and  
$f(\beta,a)$ is an arbitrary solution of the GG-system such that 
$f$ is regular on the hyperplane $\Pi_v\subset V\,$; then 
$f$ satisfy the following additional relation:
$$
\sum_{\om\in\A_v}\langle\lambda,\om\rangle a_\om
{\partial f(\beta,a)\over\partial a_{\om+v}}=0\quad
\mbox{for}\ \beta\in\Pi_v\,,
$$
where $\lambda\ne 0$ is the vector of the dual space $V^*$, 
orthogonal to $L_v\,$. 
\end{Th} 

\subsection{Example of resonance GG-system}
Let $\CC^N=\CC^p\otimes\CC^n$ be the space of $p\times n$--matrices
$a=\|a_{ij}\|\,$ ($p<n$), 
$\A$ be the set of $N=pn$ vectors $\om_{ij}=e_j+d_i\in\CC^n\oplus\CC^p$, 
where $e_1,\dots,e_n$ and $d_1,\dots,d_p$ are bases in $\CC^n$ and 
$\CC^p$ respectively, $V\subset\CC^n\oplus\CC^p$ be the linear space 
generated by the vectors $\om_{ij}\,$. 

The  GG-system in the space of functions 
$f(\alpha,\beta;\,a)$, 
$(\alpha,\beta)\in V\,$, $a\in\CC^N$, associated with the set $\A$,
has the following form:
$$
{\partial f(\alpha,\beta;\,a)\over\partial a_{ij}}
=f(\alpha-e_j,\beta-d_i;\,a)\ ,\quad i=1,\dots,p,\ j=1,\dots,n
$$
\equation\label{GGres}
\sum_{i=1}^p a_{ij}{\partial f\over\partial a_{ij}}=
  \alpha_j f\ ,\quad j=1,\dots,n
\endequation
$$
\sum_{j=1}^n a_{ij}{\partial f\over\partial a_{ij}}=
  \beta_if\ ,\quad i=1,\dots,p
$$
Note that the condition $(\alpha,\beta)\in V$ is equivalent to 
the relation $\sum\alpha_j=\sum\beta_i\,$. 

It is easily verified that all vectors $v_{ii'}=d_{i'}-d_i\in V$ 
are consistent with $\A\,$. We have: $\A_{v_{ii'}}=
\mbox{$\{\om_{ij}\bigm|j=1,\dots,n\}$}\,$, and the hyperplane 
$\Pi_{v_{ii'}}$ is defined by the equation $\beta_i=-1$ and 
therefore does not depend on $i'\,$. 

\begin{Prop}
Any GG-function $f(\alpha,\beta;\,a)$ that is regular for $\beta_i=-1\,$, 
satisfy for $\beta_i=-1$ the following additional relations:
$$
\sum_{j=1}^n a_{ij}{\partial f\over\partial a_{ij}}=0\ ,
\quad i'=1,\dots,p\,,\ i\ne i'\ .
$$ 
\end{Prop}

We denote restrictions 
of GG-functions $f(\alpha,\beta;\,a)$ 
to the plane in $V\,$, defined by the equations $\beta_i=-1\,$, 
$i=1,\dots,p\,$, by $\varphi(\alpha,a)$. 

\begin{Prop}
Functions $\varphi(\alpha,a)$ satisfy the following system of equations:
$$
\sum_{i=1}^p a_{ij}{\partial \varphi\over\partial a_{ij}}
=\alpha_j\varphi\ ,\quad j=1,\dots,n
$$
\equation\label{resphi} 
\sum_{j=1}^n a_{ij}{\partial\varphi\over\partial a_{i'j}}=-\delta_{ii'}\varphi
\endequation
$$
{\partial^2\varphi\over\partial a_{ij}\partial a_{i'j'}}=
{\partial^2\varphi\over\partial a_{ij'}\partial a_{i'j}}
$$
\end{Prop}

\noindent{\bf Remark.}\ According to \cite{Gelf} the system (\ref{resphi}) 
is the hypergeometric system of equations, associated with 
the Grassmanian $G_{p,n}$ of $p$-dimensional subspaces of $\CC^n$.
We see that this system arises as resonance case of general GG-system 
(\ref{GG1}),(\ref{GG2}).

\section{GG-functions associated with an arbitrary complex Lie group
\protect\cite{Lie}}
We give here a definition of GG-systems associated with an arbitrary 
complex Lie group $G$. This definition includes the definition of 
section 1 as a particular case. Earlier M.Kapranov \cite{K} generalized 
the definition of $\A$-hypergeometric system onto arbitrary reductive 
Lie groups. 

Let $G$ be an arbitrary connected complex Lie group and let 
$v\mapsto vg$, $v\mapsto gv$, $g\in G$ be its actions 
by the right and left translations in the space of functions on $G\,$. 
Denote by $\cal H$ the space of analytic functions $v$ on $G$ such that 
the linear space generated by the right translations $vg\,$, $g\in G\,$, 
is finite-dimensional. Evidently $\cal H$ is an algebra with respect 
to multiplication of functions, and $\cal H$ is invariant under the 
right and left translations of $G\,$. 

Let $H\subset {\cal H}$ and $V\subset{\cal H}$ be  
linear subspaces that are invariant under the right translations.
Suppose that $H$ is finite dimensional and $V$ is closed under 
multiplication by elements of $H\,$. 

\begin{Def}
The following system of equations in the space of functions 
$f(v,h)$ on $V\times H$ is called a GG-system associated with 
$V$ and $H$:
\equation\label{7.1}
f(vg,hg)=f(v,h)\quad \mbox{for every}\ g\in G\ ,
\endequation
\equation\label{7.2}
{\partial f(v,h)\over\partial h^0}=f(h^0v,h)\quad \mbox{for every}\ h^0\in H
\endequation
where $\displaystyle{\partial f(v,h)\over\partial h^0}=
{d\over dt}f(v,h+th^0)\Big|_{t=0}\,$.
A solution $f(v,h)$ of equations (\ref{7.1}),(\ref{7.2}) linear with 
respect to $v$ and analytical with respect to $h$ is called 
a GG-function associated with $V$ and $H$. 
\end{Def}

\noindent{\bf Remarks.}\ $1^{\rm o}$ If $h_1(g),\dots,h_N(g)$ is 
a fixed basis in $H$ then the equations (\ref{7.2}) have the following 
form in the coordinates $(x_1,\dots,x_N)$ corresponding to this basis: 
\equation\label{coo7.2}
{\partial f(v,x)\over\partial x_j}= f(h_jv,x)\ ,\ j=1,\dots,N
\endequation
$2^{\rm o}$ If $B$ is a basis in $V\,$, then GG-functions may be 
regarded as functions on $B\times H\,$. 
\smallskip

The GG-system (\ref{GG1}),(\ref{GG2}) is a special case of the system 
(\ref{7.1}),(\ref{7.2}) when $G$ is an additive group ${\bf C}^n$, 
$H$ is the linear space with the basis $h_i(t)=e^{\langle\om^i,t\rangle}$, 
$t\in{\bf C}^n$, and $V$ is linearly generated by all vectors 
$h_\beta(t)=e^{\langle\beta,t\rangle}$, $\beta\in ({\bf C}^*)^n$. 
\begin{Prop}
The following integrals formally satisfy the system (\ref{7.1}),(\ref{7.2}):
\equation\label{7.3}
f(v,h)=\int\limits_C e^{h(g)}v(g)\,dg\ ,
\endequation
where $dg$ is a right-invariant holomorphic volume form 
and $C\subset G$ is a cycle of real dimension equal to $\dim_{\bf C}G\,$. 
\end{Prop}

Let us consider a case when the orbits of $G$ on the space $H$ are 
algebraic submanifolds. Let $P\subset{\bf C}[\xi_1,\dots,\xi_N]$ 
be the subring of all polynomials $p(\xi)=p(\xi_1,\dots,\xi_N)$ 
such that $p\bigl(h_1(g),\dots,h_N(g)\bigr)\equiv 0\,$, 
where $h_1,\dots,h_N$ is a fixed basis in $H\,$. 

Fix an arbitrary finite-dimensional subspace $L\subset {\cal H}$ 
and denote its dual by $L^*$. 
\begin{Def}
The following system of equations in the space of functions 
$F(h)$ on $H$  with values in $L^*$ is called the general GGZ-system 
associated with $L$ and $H\,$:
\equation\label{7.4}
\langle F(hg),v\rangle=\langle F(h),vg^{-1}\rangle\quad
\mbox{for every}\ v\in L\ \mbox{and}\ g\in G
\endequation
\equation\label{7.5}
p({\partial\over\partial x})F=0\quad\mbox{for every}\ p\in P\ .
\endequation
Solutions of (\ref{7.4}),(\ref{7.5}) in the class of analytic functions 
on $H$ are called GGZ-functions associated with $L$ and $H\,$. 
\end{Def}

The $\A$-hypergeometric system (\ref{HS1}),(\ref{HS2}) is a special 
case of (\ref{7.4}),(\ref{7.5}) when $G$ is an additive group ${\bf C}^n$, 
$H$ is the linear space with the basis $h_i(t)=e^{\langle\om^i,t\rangle}$, 
$t\in{\bf C}^n$, $i=1,\dots,N$, and $V$ is the one-dimensional space 
generated by $e^{\langle\beta,t\rangle}$. 

Let $f(v,h)$ be a function on $V\times H$ linear with respect to $v$, 
and let $L\subset V\,$. Define a function $F(h)$ on $H$ with values 
in $L^*$ by the equality:
$$
\langle F(h),v\rangle =f(v,h)\quad \mbox{for any}\ v\in L
$$
\begin{Prop}
If $f(v,h)$ is a GG-function associated with $V$ and $H\,$, then 
$F(h)$ is a GGZ-function associated with $L$ and $H\,$. 
\end{Prop}

\noindent{\bf Remark.}\ It is natural to define the q-analog 
of GG-system (\ref{7.1}),(\ref{7.2}) as the system consisting of 
equations (\ref{7.1}) and of the equations obtained from equations (\ref{coo7.2}) 
by replacing the operators $\partial/\partial x_j$ by corresponding 
q-differential operators. This system depends not only on the subspaces 
$V$ and $H$ but also on a choice of the basis $h_1,\dots,h_N$ in $H\,$. 
Similarly one can define the q-analog of the general GGZ-system. 
The following functions formally satisfy the q-analog of the GG-system: 
$$
f(v,h)=\int\limits_C\prod_{j=1}^N \exp_q\bigl(x_jh_j(g)\bigr)v(g)\,dg\ ,
$$
where $\exp_q$ is the q-exponential and $dg$ and $C$  are defined 
in the same way as in (\ref{7.3}).


\begin{thebibliography}{44}
\bibitem{Gelf} Gelfand I.M. General theory of hypergeometric functions.
  Soviet Math.~Doklady, vol.~33 (1986), 573--577 
\bibitem{G1} Gelfand I.M., Gelfand S.I. Generalized hypergeometric systems.
  Soviet Math.~Doklady, vol.~33 (1986), 279--283
\bibitem{G2} Gelfand I.M., Graev M.I. Duality theorem for general 
hypergeometric functions.
  Soviet Math.~Doklady, vol.~34 (1987), 9--13
\bibitem{G36} Gelfand I.M., Graev M.I. Hypergeometric functions 
associated with the Grassmanian $G_{3,6}$.
  Soviet Math.~Doklady, vol.~35 (1987), 298--303
\bibitem{G36a} Gelfand I.M., Graev M.I. Hypergeometric functions 
associated with the Grassmanian $G_{3,6}$.
  Soviet Math.~Sbornik, vol.~180 (1989), 3--38
\bibitem{GVZ}  Vassiliev V.A., Gelfand I.M., Zelevinsky A.V. 
General hypergeometric functions on complex Grassmanians.
  Funct.~Anal.~Appl., vol.~21 (1987), 19--31 
\bibitem{GGZ} Gelfand I.M., Graev M.I., Zelevinsky A.V. 
Holonomic systems of equations and series of hypergeometric type.
  Soviet Math.~Doklady, vol.~36 (1988), 5--10
\bibitem{GGf1} Gelfand I.M., Graev M.I. GG-functions of one variable. 
  Russian Acad.~Sci.~Dokl.~Math., vol.~55 (1997)
\bibitem{GGf2} Gelfand I.M., Graev M.I. GG-functions of multiple 
variables and their connections with general hypergeometric functions.
  Russian Acad.~Sci.~Dokl.~Math., vol.~55 (1997)
\bibitem{GGf3} Gelfand I.M., Graev M.I. GG-functions of multiple 
variables and their connections with general hypergeometric functions.
  Russian Math.~Surveys, vol.~52, no.~4 (1997), 639--684 
\bibitem{Loop2} Pham. Formules de Picard-Lefschetz generalisees 
et ramification des integrales. Bull.~Soc.~Math.~France {\bf 93} (1965), 
333--367
\bibitem{AI} Aomoto K., Iguchi K. On quasi hypergeometric functions. 
To appear in Methods and Applications of Analysis. 
\bibitem{K} Kapranov M. Integrable systems and algebraic geometry.
In: Proc.~of the Tauguchi symp. 1997. M.-N.Saito, I.Shimazu, K.Ueno eds.
 World Scientific, 1998. P.236--281.
\bibitem{Lie} Gelfand I.M., Graev M.I. The special functions associated 
with complex Lie groups. Doklady Rossiyskoy Akademii Nauk, v.~364, No.~2 
(1999), 151--154 (in Russian). English Transl.: Russian Acad.~Sci.~Dokl.~Math., 
v.~59, No.~1 (1999)
\end{thebibliography}
\end{document}